\documentclass[11pt]{amsart}
\usepackage[T1]{fontenc}
\usepackage{amssymb,amsmath,amsthm}

\newtheorem{theorem}{Theorem}[subsection]

\newtheorem{claim}[theorem]{Claim}

\newtheorem{corollary}[theorem]{Corollary}

\newtheorem{definition}{Definition}[section]

\newtheorem{lemma}[theorem]{Lemma}

\newtheorem{proposition}[theorem]{Proposition}
\newtheorem{remark}[theorem]{Remark}

\newenvironment{postulate}[1][Consistent Amalgamation.]{\begin{trivlist}
\item[\hskip \labelsep {\bfseries #1}]}{\end{trivlist}}

\begin{document}

\def\Ind#1#2{#1\setbox0=\hbox{$#1x$}\kern\wd0\hbox to 0pt{\hss$#1\mid$\hss}
\lower.9\ht0\hbox to 0pt{\hss$#1\smile$\hss}\kern\wd0}
\def\ind{\mathop{\mathpalette\Ind{}}}
\def\Notind#1#2{#1\setbox0=\hbox{$#1x$}\kern\wd0\hbox to 0pt{\mathchardef
\nn=12854\hss$#1\nn$\kern1.4\wd0\hss}\hbox to
0pt{\hss$#1\mid$\hss}\lower.9\ht0 \hbox to
0pt{\hss$#1\smile$\hss}\kern\wd0}
\def\nind{\mathop{\mathpalette\Notind{}}}

\def\thind{\mathop{\mathpalette\Ind{}}^{\text{\th}} }
\def\nthind{\mathop{\mathpalette\Notind{}}^{\text{\th}} }
\def\uth{\text{U}^{\text{\th}} }

\def\spec{\mathop{Spec}}

\def\field{^{field}}
\def\qf{^{qf}}


\title{th-Forking, Algebraic Independence and Examples of Rosy Theories.}
\author{Alf Onshuus}


\maketitle

\section{Introduction}
In \cite{Onshuus} we developed the notions of \th-independence and
\th-ranks which define a geometric independence relation in a
class of theories which we called ``rosy''. We proved that rosy
theories include simple and o-minimal theories and that for any
theory for which the stable forking conjecture was true,
\th-forking coincides with forking independence.

In this article, we continue to study properties of \th-forking
and examples of rosy theories. In the first section we recall the
definitions and basic results proved in \cite{Onshuus}. In section
2 we study alternative ways to characterize rosy theories and we
prove the coordinatization theorem, which will prove very useful
when studying examples of superrosy theories. 
In section 3 we prove that given any theory $T$
with weak elimination of imaginaries and in which algebraic
independence defines a geometric independence relation, $T$ is
rosy and \th-forking agrees with algebraic independence.
Finally, in section 4 we
study two examples of rosy theories. We 
prove that pseudo real closed fields (PRC-fields) are rosy and that 
whenever one has a large differential field whose restriction to the
language of rings is a model complete rosy field, 
the model companion (as defined
by Tressl in \cite{tressl}) is also rosy. This is an extension of
the result we had in \cite{phd} for closed ordered differential
fields. 

Throughout this article we assume familiarity with the
terminology, definitions and basic results in stability and
simplicity theory. We work with first order theories which will
generally be denoted by $T$. As is usual in stability theory,
unless otherwise specified, we work inside a monster model
$\mathcal{C}^{eq}$ of $T$. By convention lower case letters
$a,b,c,d$ will in general represent tuples (of imaginaries) unless
otherwise specified, and upper case letters will represent sets.
Greek letters such as $\delta, \sigma, \psi, \phi$ will be used
for formulas. Given a set $A$ and a tuple $a$, the type of $a$
over $A$ will be denoted by $tp(a/A)$) and we will abbreviate
``$tp(a/Ab)$ does not fork over $A$'' as $a\ind _A b$.

\subsection{Preliminaries and notation}

We recall the definitions of \th-forking and rosy theories and the
main results proved in \cite{Onshuus}. The main definitions we
work with are the following:

\begin{definition}

A formula $\delta (x,a)$ \emph{strongly divides over $A$} if
$tp(a/A)$ is non-algebraic and $\{ \delta (x,a') \}_{a'\models
tp(a/A)}$ is $k$-inconsistent for some $k\in \mathbb{N}$.

We say that $\delta (x,a)$ \emph{\th-divides over $A$} if we can
find some tuple $c$ such that $\delta (x,a)$ strongly divides over
$Ac$.

A formula \emph{\th-forks over $A$} if it implies a (finite)
disjunction of formulae which \th -divide over $A$.

As is done with the standard forking, we say that the type $p(x)$
\emph{\th-divides over $A$} if there is a formula in $p(x)$ which
\th-divides over $A$; \th-forking is similarly denoted. We say
that $a$ is \emph{\th-independent of $b$ over $A$}, denoted
$a\thind_A b$, if $tp\left( a/Ab\right) $ does not \th -fork over
$A$.
\end{definition}

As we mentioned before, \th-forking defines an independence
relation in a large class of theories called rosy theories which
includes simple and o-minimal structures. Before we are able to
give the definition of rosy theories we must first define the
class of ranks that is associated with \th-forking.

\begin{definition}
Given a formula $\phi ,$ a set $ \Delta $ of formulas in the
variables $x;y$, a set of formulae $\Pi$ in the variables $y;z$
(with $z$ possibly of infinite length) and a number $k,$ we define
\th $ \left( \phi ,\Delta, \Pi ,k\right) $ inductively as follows:
\begin{enumerate}
\item  \th $\left( \phi ,\Delta , \Pi, k\right) \geq 0$ if $\phi $ is
consistent.

\item  For $\lambda $ limit ordinal,
\th $ \left( \phi ,\Delta, \Pi ,k\right) \geq \lambda $ if and
only if \th $ \left( \phi ,\Delta , \Pi, k\right) \geq \alpha $
for all $\alpha <\lambda $

\item  \th $ \left( \phi ,\Delta, \Pi ,k\right) \geq \alpha +1$ if and
only if there is a $ \delta \in \Delta $, some $\pi(y;z) \in \Pi$
and parameters $c$ such that

\begin{enumerate}
\item  \th $ \left( \phi \wedge \delta \left( x,a\right) ,\Delta ,\Pi
,k\right) \geq \alpha $ for infinitely many $a\models \pi(y;c) $

\item  $\left\{ \delta \left( x,a\right) \right\} _{a \models \pi(y;c)}$ is
$k-$inconsistent
\end{enumerate}
\end{enumerate}

\end{definition}

A theory is called ``rosy'' if all of its \th-ranks are defined.
It is clear from the definitions that all simple theories are rosy
(for any \th-rank we can easily find a D-rank such that the value
is bigger for all formulas). We also have the following theorem.

\begin{theorem}
In a rosy theory \th-forking has all the properties of a geometric
independence relation.
\end{theorem}

We should point out that the independence theorem does not always
work for \th-forking (as evidenced by the fact that it does define
an independence relation in o-minimal theories) and we therefore
cannot conclude that \th-forking is the same as forking when
restricted to simple theories. We do know that in all stable and
supersimple theories (and in fact in any simple theory for which
\th-forking satisfies the stable forking conjecture) forking is
the same as \th-forking \footnote{Clifton Ealy has extended this
result to all simple theories that satisfy elimination of
hyperimaginaries}.

\section{Characterization of Rosy Theories and Coordinatization Theory}

In this section we study alternative ways of characterizing
\th-forking. We prove that symmetry and local character of
\th-forking are both characteristics that imply rosiness. We also
prove the coordinatization theorem (proved by Hart, Kim and Pillay
-~\cite{Hart/Kim/Pillay}- for simple theories).

\begin{remark}
Given any type $p$ over some set $B\supset A$, if $p$ is finitely
satisfiable over $A$ then it does not \th-fork over $A$.
\end{remark}

\begin{proof}
We know that if a type is finitely satisfiable it does not fork
over $A$. \emph{A fortiori}, it cannot \th-fork.
\end{proof}

For the proofs in this section we use repeatedly the following
result. The proof is very similar to its simple theoretic analogue
(\cite{Wagner} 2.3.7).

\begin{lemma}
\label{2.3.7} Let \th$(x=x, \phi, \pi, k)\geq n$ for all $n<
\omega$. Then for every linearly ordered index set $I$, there are
\begin{enumerate}
\item An indiscernible sequence $(b_ia_i:i\in I)$ such that
$\models \phi(b_i,a_j)$ for any $j\leq i$ and $\phi(x,a_i)$
\th-forks over $\{ a_j \mid j<i \}$.

\item A tuple $b$ and a $b$-indiscernible infinite sequence
$\left< a_i \right>$ such that $\models \phi(b, a_i)$ and
$\phi(x,a_i)$ \th-forks over $\{ a_j \mid j<i \}$.
\end{enumerate}
\end{lemma}

\begin{proof}
\noindent \emph{1.} By definition of \th-rank, given $p(x)$ a
partial type over $A$, if \th$(p(x), \phi, \pi, k)\geq n+1$ we can
find tuples $a,b,c$ such that $\models \phi(b,a)\wedge \pi(a,c)$,
$tp(a/Ac)$ is non algebraic, \th$(p(x)\cup \{\phi (x,a) \}, \phi,
\pi, k)\geq n$ and $\{ \phi(x,a') \}_{a'\models \pi (y,c) } $ is
$k$-inconsistent. Thus, for every $n$, we can find by induction
tuples $\left< a_i^n, b_i^n, c_i^n \right>_{i\leq n}$ such that
$b_i^n\models \phi (x,a_j^n)$ for any $j\leq i$, $\{ \phi(x,a')
\}_{a'\models \pi (y,c_i^n) } $ is $k$-inconsistent,

\noindent $tp(a_i^n/Ac_i^na_0^n, a_1^n, \dots a_{i-1}^n)$ is
non-algebraic and
\[ \text{\th} \left( \bigcup_{j\leq i} \{\phi (x,a_j^n) \}, \phi, \pi, k \right) \geq
n-i.\]

By compactness we can build a similar sequence of any given length
and using Ramsey's theorem and compactness we can find one such
sequence such that $\left< a_ib_i \right>$ is indiscernible.

\noindent \emph{2.} Is a particular case of the sequence we get
above with $b=b_{\omega}$.

\end{proof}

\begin{remark}
\label{remark2.3.7} We are actually proving something slightly
stronger. If in the previous lemma we assume that for some type
$p(x)$, \th$(p(x),\phi, \pi, k)> n$ for all $n < \omega$, then we
can actually get all the $b$'s and $b_i$'s in the conclusions
satisfying $p(x)$.
\end{remark}

\subsection{Characterizations by Local Character and Symmetry}

\begin{theorem}
\label{local character} A theory is rosy if and only if
\th-forking satisfies local character.
\end{theorem}

\begin{proof}
To prove the left to right implication, assume that all the
\th-ranks are finite. Fixing a type $p$, for any formulas $\phi,
\pi$ and any integer $k$, there must be some finite tuple $a_{\phi
\pi k}$ such that \th$(p, \phi, \pi, k)$=\th$(p\upharpoonright
a_{\phi \pi k}, \phi, \pi, k)$. Then, if we let $A$ be the union
of all such tuples, $|A|\leq |T|$ and $p$ does not \th-fork over
$A$.

For the other direction, assume that $D(x=x, \phi, \pi, k)\nless
\omega$. By \ref{2.3.7} there is some $b$ and some
$b$-indiscernible sequence $\left< a_i \right>_{i\in |T|^+}$ such
that for any $i\in |T|^+$, $\phi (x,a_i) $ \th-divides over $\{a_j
\mid j< i \}$ and $b\models \phi(x,a_i)$. Let $A=\{a_i \mid i\in
|T|^+ \}$. Now, given any $A_0\subset A$ of cardinality $|T|$, let
$\lambda $ be some ordinal in $|T|^+$ such that $A_0\subset \{a_j
\mid j< \lambda \}$. Then $\phi(b,a_\lambda)\in tp(b/A)$ and it
\th-forks over $A_0$, contradicting local character for $tp(b/A)$.
\end{proof}

\begin{theorem}
A theory is rosy if and only if \th-forking satisfies symmetry.
\end{theorem}

\begin{proof}
The left to right implication follows from \th-forking being a
geometric independence relation in rosy theories. For the other
implication, let us assume again that $D(x=x, \phi, \pi, k)\nless
\omega)$ and by \ref{2.3.7} there is some $b$, some
$b$-indiscernible sequence $\left< a_i \right>_{i\leq \omega}$ and
some $\phi (x,y) $ such that for all $i$ we have that $b\models
\phi(x,a_i)$ and $\phi(x,a_i)$ \th-divides (and therefore it
\th-forks) over $\{a_j \mid j< i \}$. By indiscernibility,
$tp(a_\omega/ba_1,\dots ,a_n,\dots )$ is finitely satisfiable over
$a_1, \dots ,a_n \dots$ and thus it does not \th-fork over $a_1,
\dots ,a_n \dots$. On the other hand $b\models \phi(x,a_\omega)$
so in any non-rosy theory, \th-forking is not symmetric.
\end{proof}

\subsection{Coordinatization}
All the definitions and results in this subsection analogues of
simple theoretic results proved in ~\cite{Hart/Kim/Pillay}; the
proofs are very similar.

\begin{definition}
\label{coordinatization} We say that a theory $T$ is finitely
coordinatized by a set of types $\mathcal{P}$ if $\mathcal{P}$ is
closed under automorphisms and for any type $p(x)$ over a tuple
$a$ there is some $n\in \omega$ and a sequence $a_0, \dots, a_n$
with $a_0=\emptyset$, $a_n=a$ and $tp(a_i/a_{i-1})\in
\mathcal{P}$.
\end{definition}

\begin{definition}
Given any type $p\in S(T)$, we say that $p$ is \emph{rosy} if all
of its  extensions satisfy local character for  \th-forking.
\end{definition}

\begin{lemma}
\label{local local character} Let $p(x)$ be any type. Then
\th$(p(x),\delta, \pi, k)< \omega$ for all formulas $\delta$,
$\pi$ and all $k\in \mathbb{N}$ if and only if $p(x)$ is rosy.
\end{lemma}

\begin{proof}
Using \ref{remark2.3.7} the result follows from the proof of
\ref{local character}.
\end{proof}

\begin{theorem}
\label{2.5.2.2} If $T$ is coordinatized by rosy types, then it is
rosy.
\end{theorem}

\begin{proof}
Let $a$ be any tuple, and $A$ be a set; we can assume $a$ itself
is a coordinatizing sequence (after maybe adding entries to $a$)
of size, say, $n$. We will prove that for any $\phi, \pi, k$,
\th$(tp(a/A), \phi, \pi, k)$ is finite by induction on $n$. If
$n=1$ the result follows from \ref{local local character}. So let
us assume $a=bc$ where $b$ is a (coordinatizing) sequence of size
$n-1$ and $c$ is a singleton. By induction hypothesis, there is
some $A_0\subset A$ of size less than $|T|$ such that $tp(b/A)$
does not \th-fork over $A_0$. By local character (maybe after
increasing $A_0$ but keeping its size smaller than $|T|$) we can
assume $tp(c/Ab)$ does not \th-fork over $A_0b$. By partial left
transitivity (\cite{Onshuus} lemma 2.1.5) we get that $tp(bc/A)$
does not \th-fork over $A_0$.
\end{proof}

\begin{definition}
A theory $T$ is \emph{ \th-minimal } if for any model $M\models
T$, any $a\in M^1$ and any $A\subset M$, $tp(a/A)$ \th-forks over
$\emptyset$ if and only if $a$ is algebraic over $A$.
\end{definition}

\begin{corollary}
\label{cor-to-coord} Given a complete theory $T$, let $M\models T$
be a model of $T$. Then, if for all elements $a\in M^1$ and all
sets $A\subset M$ $\uth (a/A)\leq 1$, $T$ is superrosy. In
particular, any \th-minimal theory is superrosy.
\end{corollary}

\begin{proof}
Clearly the set of all types of single elements coordinatizes the
theory. Being \th-minimal is equivalent to have the $\uth$-rank of
any such type (and therefore all the \th-ranks) be at most 1 and
by \ref{local local character} any type in $S_1(T)$ is rosy. By
theorem \ref{2.5.2.2} $T$ is rosy and by Lascar's inequalities for
rosy theories it is in fact superrosy.

Finally, let $T$ be \th-minimal, $p\in S_1(T)$ be a type over $A$,
$q$ be a type over $B\supset A$ extending $p$ and let $a$ be an
element satisfying $q$. If $q$ is not algebraic, then $q$ does not
\th-fork over the empty set and it is rosy. If $q$ is algebraic,
then there is some finite $C\subset B$ such that $tp(a/C)$ is
algebraic and $q$ does not \th-fork over $C$.
\end{proof}

\section{\th-Forking and Algebraic Independence in Non-Saturated Models}

Besides simple theories, most examples of
theories with an independence relation are those for which algebraic
independence satisfies the Steinitz exchange property.
In this section we prove that if algebraic independence is a
geometric independence relation and one has elimination of
imaginaries the theory is rosy and algebraic independence
corresponds to \th-forking.

This results were proved in \cite{phd} when we were studying what
(weakened) version of amalgamation of independent types one could
prove for \th-forking in a general rosy theory. It is known (see
for example \cite{Kim}) that the only independence relation
satisfying symmetry, local character, transitivity and having
amalgamation for independent types (the independence theorem) was
the non forking relation inside a simple theory. This means that
the only rosy theories in which \th-forking can have the
independence theorem are simple theories in which \th-forking
agrees with forking.

However, in our main non-simple examples (o-minimal theories) the
following amalgamation theorem is true:
\begin{postulate}
\label{o-min-am}
Let $p(x)$ be a complete type over some set $A$. Let $p(x,a)$ and
$q(x,b)$ be two non-\th-forking extensions of $p(x)$ such that
$a\thind _A b$. Then, either $p(x,a)\cup p(x,b)$ is inconsistent,
or it is a non-\th-forking extension of $p(x)$.
\end{postulate}

All the examples of rosy theories mentioned here have consistent
amalgamation. The results in this section were used to find an
example of a rosy theory that does not have consistent
amalgamation. The details of the construction can be found in
\cite{phd}.

\subsection{\th-forking in non-saturated models}

Up until now we have followed the simple theoretic approach of
working inside a large saturated model $\mathcal{C}\models T$. In
simple theories this is necessary to understand the behavior of
forking-independence in $T$, mainly because the definition of
forking requires indiscernible sequences and the existence of such
sequences varies inside different models of $T$. We prove that,
even though we still need some saturation to fully understand the
behavior of \th-forking, the amount of saturation we need is far
from what we need to study simple theories.

\begin{definition}
Let $T$ be a theory, $M$ be a model of $T$ and $\mathcal{C}$ be a monster
model of $T$ containing $M$.

Given a set $A$ in $M$, the \emph{algebraic closure of $A$ in
$M$}, which we denote as $acl^M(A)$, is the set of elements $b$ in
$M$ such that $tp(b/A)$ has infinitely many realizations in $M$.

A model $M\models T$ is \emph{weakly $\omega$-saturated} if for any $c\in M$
any formula of the form
\[ \phi \left( x_1,x_2, \dots, x_n; c \right) \] if there are
$a_1, a_2, \dots , a_n$ in $\mathcal{C}$ such that
$\mathcal{C}\models \phi(a_1, \dots, a_n;c)$ and for any $i$ $a_i\notin
acl(a_{i+1}, \dots a_n,c)$ then there are $b_1, b_2, \dots ,b_n$ in $M$ such
that $M\models \phi(b_1, \dots, b_n;c)$ and for any $i$ $b_i\notin
acl^M(b_{i+1}, \dots b_n,c)$.
\end{definition}

\begin{theorem}
Let $M$ be a model of $T$ that is weakly $\omega$-saturated. Then
 \th-forking dependence (between tuples in $M$) is allways witnessed by
elements in $M$.

This is, if $a,b\in M$, $A\subset M$ and $a\nthind _A b$ then
there are formulas
$\psi_i (x,b_i)$, $b_i\in M$ and tuples $c_i\in M$ such that
\begin{itemize}
\item $tp(a/Ab)\models \bigvee_{i=1}^m \psi_i(x,b_i)$,
\item $\{ \psi_i(x,b_i') \}_{b_i'\models tp(b_i/Ac_i)}$ is $k$-inconsistent and
\item $tp(b_i/Ac_i)$ has infinitely many realizations in $M$.
\end{itemize}
\end{theorem}

\begin{proof}
$k$-inconsistency can be witnessed by a formula so all three of the statements
 follow from the definition of \th-forking and weak $\omega$-saturation.
\end{proof}

\begin{corollary}
\label{nonsaturated}
A theory $T$ is \th-minimal if and only if for any weakly $\omega$-stable
model $M$, any tuples $a,b\in M^1$ and any finite $A\subset M$, $a\nthind _A b$
if and only if $tp(a/Ab)$ is algebraic in $M$ and $tp(a/A)$ has infinitely
many realizations (in $M$).
\end{corollary}

\subsection{\th-Forking and algebraic independence}

\begin{definition}
Following Pillay, we say that a theory $T$ has
\emph{geometric elimination of imaginaries} if for any model $M\models
T$ and any imaginary element $e\in dcl^{eq}(M)$, there is some
finite $x\subset M$ such that $acl^{eq}(e)=acl^{eq}(x)$. We say
that a model $M$ weakly eliminates imaginaries if $Th(M)$ does.
\end{definition}

\begin{theorem}
\label{uth1}
Let $M$ be a weakly $\omega$-stable model which satisfies the Steinitz
exchange property
for the algebraic closure and which has weak elimination of
imaginaries. Then $Th(M)$ is \th-minimal and has global
$\uth$-rank 1. Moreover,   $\uth(A)=\dim^{\rm alg}(A)$.
\end{theorem}

\begin{proof}

By \ref{cor-to-coord} and \ref{nonsaturated} it is enough to show that for
any $a\in M^1$ and any $A\subset M$, $\uth(a/A)\leq 1$ .

We proceed by contradiction. Suppose we have some $b\in M^1$ such
that $\uth(b)\geq 2$. By definition there are finite $A,C\subset
M$ such that $p(x,A):=tp(b/A)$ is non-algebraic and there is some
formula $\phi(x,A)\in p(x,A)$ such that
\[ \{ \phi(x,A') \}_{A'\models tp(A/C)}\]  is $k$-inconsistent and
$tp(A/C)$ is non-algebraic. Using weak elimination of imaginaries,
we may replace $A$ and $C$ by tuples of
elements in $M^1$.

\begin{claim}
We may assume without loss of
generality that $b$ is not algebraic over $A\cup C$.
\end{claim}
\begin{proof}
Since $\phi(\mathcal{C},A)$ is non algebraic it is unbounded so
there is some $b'\in \mathcal{C}$ satisfying $\phi(x,A)$ and such that
$b'$ is not algebraic over $A\cup C$. By weak $\omega$-saturation, we can find
such a $b$ in $M$.
\end{proof}

$\{ \phi(x,A') \}_{A'\models tp(A/C)}$ is $k$-inconsistent so
$tp(A/Cb)$ is algebraic. Since $tp(A/C)$ is not algebraic by
definition and $A$ and $C$ are tuples, we can apply Steinitz
exchange property repeatedly and eventually we get that $tp(b/AC)$
is algebraic, a contradiction.
\end{proof}

\section{ Two Examples of Rosy Theories}

\subsection{Pseudo Real Closed Fields}

In \cite{cherlin/van/mc}, Cherlin, van den Dries and Macintyre
proved that the theory of a bounded pseudo-algebraically closed
field\footnote{A field $F$ is pseudo-algebraically closed (PAC) if
and only if any absolutely irreducible variety defined by
polynomials with coefficients in $F$ has an $F$-rational point. A
PAC field is bounded if the absolute Galois group is small: this
is, if for any $n\in \mathbb{N}$ the number of elements in
$Gal(acl(F):F)$ of degree $n$ is finite.} is decidable and
therefore behaves nicely in a model theoretic sense. In a series
of papers (see \cite{hrushovskimanuscript}, \cite{chatzidakis} and
\cite{chatzidakis/pillay}), Hrushovski, Chatzidakis and Pillay
proved that any such theory (the theory of a bounded
pseudo-algebraically closed field) is simple; even more, they
proved that a pseudo-algebraically closed field is simple if and
only if it is bounded.

In \cite{Prestel}, Prestel defined the theory of pseudo real
closed fields (PRC-fields) which includes ordered (or orderable)
fields that are ``close'' to being real closed fields (in the same
sense in which PAC fields are ``close'' to being algebraically
closed). Prestel's definition states that a field $F$ is a
PRC-field if any absolutely irreducible variety that has a
$K$-rational point in every real closed field $K$ containing $F$,
has an $F$-rational point. This definition is a strict
generalization of PAC-fields (there are no real closed fields
containing a PAC-field so the condition of having a rational point
in every real closure becomes vacuous) and works nicely when
trying to describe the possible absolute Galois groups of a
PRC-field.

We work with a definition closer to the one used in
\cite{hrushovskimanuscript} for PAC-fields. We work inside a two
sorted structure; let $(D,F)$ be a pair of ordered fields such
that $D$ is the real closure of $F$ and such that if $V$ is an
open cell in $D^n$ such that the (topological) closure of $V$ is
defined as the ``real'' zeros of an absolutely irreducible variety
definable over $F$, then $V$ has $F$-rational points. Let
$T=Th((D,F))$.

The first thing we show is that this definition is not far from
the one stated in \cite{Prestel}.

\begin{claim} For any PRC-field (in the sense of Prestel)
$F$, if $D$ is the real closure of $F$, then $(D,F)\models T$.
\end{claim}

\begin{proof}

In \cite{Schmid} Theorem 2.1 Schmid proved that if $F$ is a
PRC-field and $V$ is an absolutely irreducible affine variety
definable over $F$ then the set of $F$-rational points in $V$ is a
dense subset (in the order topology) of the set of $D$-rational
points in $V$. The claim follows.
\end{proof}

In the rest of the section we prove that if $F$ is a bounded
  PRC field,
then $T$ is rosy. In fact we prove the following theorem.

\begin{theorem}
\label{minimalityPRC} Let $F$ be a bounded,   PRC field. Let $T^*$
be $Th((F,a)_{a\in N})$ where $N$ is an elementary submodel of
$F$. Then $T^*$ has elimination of imaginaries and algebraic
closure satisfies Steinitz exchange property.
\end{theorem}

Both the proofs of elimination of imaginaries and Steinitz
exchange property rely on the proofs given by Hrushovski for
perfect PAC fields in \cite{hrushovskimanuscript}.

\begin{lemma}
\label{PRCPAC} Let $F$ be a PRC field, let $D$ be its real closure
and let $i=\sqrt{-1}$. Then $F(i)$ is a perfect PAC-field and
$D(i)$ is its algebraic closure.
\end{lemma}

\begin{proof}
Let $V$ be an absolutely irreducible variety definable over
$F(i)$. We want to show that $V$ has $F(i)$-rational points.

Being irreducible is a property that can be determined by looking
at the restrictions of $V$ inside an open affine cover, so we can
assume that $V$ is an affine variety, $V=\spec(F(i)[\bar{x}]/I)$
where $I$ is an ideal generated by polynomials $\{p_1,\dots p_m
\}$.

For any $a\in F(i)$ let $a^c$ be the complex conjugate of $a$ and
let $a^r,a^{im}$ be elements in $F$ such that $a=a^r+ia^{im}$. We
define $p_i^c$, $V^c$ and $V^R$ in the following way. Let
$\bar{a}:=\left< a_i\right> _{i\leq n}$ be a tuple in $F(i)$. Let
$V^c$ be the variety  such that for all $a$, $\bar{a}\in V$ if and
only if $\overline{a^c}:=\left< a^c_i\right> _{i\leq n}$ is a
tuple in $V^c$; let $p_i^c$ be a polynomial such that for all $x$,
$p_i(x)=0$ if and only if $p_i^c(x^c)=0$ so that
\[V^c=\spec\left( F(i)[\bar{x}]/\left< p_1^c,\dots p_m^c\right>  \right).\]
Finally, let $V^R$ be the smallest variety such that $ \bar{a}\in
V$ if and only if
\[ \overline{a^R}:=\left< a_i^r,a_i^{im}\right> _{i\leq n}\in V^R.\] In particular,
if a tuple $\left< a_j,b_j\right> $ is in $V^R$ then for any
$q(x)\in I$, $q(\left< a_j\right> +i\left< b_j\right> )=0$ and
$q^c(\left< a_j\right> -i\left< b_j\right> )=0$.

\begin{claim}
\label{isomorphism} The map
\[ \sigma: V\times V^c \rightarrow F(i)^{2n} \] defined by sending
\[\left( \left< a_j\right> , \left< a^c_j\right>  \right) \text{ to } \left( \frac{\left< a_j\right> +\left< a_j^c\right> }{2}, \frac{\left< a_j\right> -\left< a_j^c\right> }{2i} \right)=\left( \left< a_j^r,a_j^{im}\right> _{j\leq n}\right) \] is
an isomorphism from $V\times V^c$ onto $V^R$.
\end{claim}

\begin{proof}
The composition of the natural map $V\rightarrow V\times V^c$ with
$\sigma$ sends a tuple $\bar{a}$ in $V$ to $\overline{a^R}$. By
definition the image of $\sigma$ contains $V^R$.

Now, given any tuple $\left< a_j,b_j\right> _{i\leq n}$ in $V^R$,
consider the map $\tau$ that sends $\left< a_j,b_j\right> $ to the
pair $(\left< a_j+ib_j\right> ,\left< a_j-ib_j\right> )$. For any
$q(x)\in I$, $q(a_j+ib_j)=0$ and $q(a_j-ib_j)=0$; by definition,
$(\left< a_j+ib_j\right> ,\left< a_j-ib_j\right> ) \in V\times
V^c$ and taking the composition of the two maps we find that
$\tau$ is the inverse of $\sigma$.
\end{proof}

\begin{corollary}
If $V$ is irreducible, so is $V^R$.
\end{corollary}
\begin{proof}
$V^c$ is isomorphic to $V$ and therefore irreducible. The product
of two irreducible varieties is irreducible.
\end{proof}

To finish the proof of the lemma, we need to show that $V^R$ is
definable over $F$.

\begin{claim}
\label{conjugates} For any $a,b$, $\frac{b^c+a^c}{2}=\left(
\frac{a+b}{2} \right)^c$ and $\frac{b^c-a^c}{2i}=\left(
\frac{a-b}{2i}\right) ^c$.
\end{claim}

\begin{proof}
The first equality in the claim is clear. For the second part, a
simple calculation shows that for any $d\in D(i)$,
$(d/i)^c=(-d^c)/i$; taking $d=a-b$ proves the claim.
\end{proof}

\begin{corollary}
\label{conjugatestuples} Given any $\bar{a}$ and $\bar{b}$, if
$\left( \bar{a},\bar{b}\right)$ is in $V^R$ then the conjugate
$\left( \overline{a^c},\overline{b^c}\right)$ is in $V^R$ and
therefore $V^R$ is definable over $F$.
\end{corollary}

\begin{proof}
Since addition and multiplication by constants can be done in each
coordinate of a tuple, we can apply claim \ref{conjugates} to
$\bar{a}$ and $\bar{b}$. By \ref{isomorphism} there are $(\bar{c},
\bar{d})$ in $V\times V^c$ such that
\[\left( \bar{a},\bar{b}\right) = \left( \frac{\bar{c}+\bar{d}}{2}, \frac{\bar{c}-\bar{d}}{2i}\right) .\]
By definition, $(\overline{d^c}, \overline{c^c})$ is also in
$V\times V^c$ and using \ref{isomorphism} once again,
\[ \left( \frac{\overline{d^c}+\overline{c^d}}{2}, \frac{\overline{d^c}-\overline{c^c}}{2i}\right) =\left( \overline{a^c},\overline{b^c}\right)\] is in
$V^R$. By hypothesis $V^R\cong V\times V^c$ is definable over
$F(i)$ so $V^R$ is in fact definable over $F$. \end{proof}

To finish the proof of lemma \ref{PRCPAC} note that $D(i)$ is
algebraically closed so $V$ has $D(i)$-rational points. Taking the
image of such points in $V^R$ (the real and imaginary components
of the $D(i)$-rational points contained in $V$) we get that
$V^R(D)$ is non-empty. By definition of PRC-fields there is some
$F$-rational point $\left< a_j,b_j\right> $ in $V^R$ and $\left<
a_j+ib_j\right> $ is an $F(i)$-rational point in $V$.

Since $V$ was any irreducible $F(i)$-definable variety, every
absolutely irreducible variety definable over $F(i)$ has
$F(i)$-rational points. \end{proof}

To prove theorem \ref{minimalityPRC} we follow the proofs of weak
elimination of imaginaries for PAC-fields (corollary 3.2 in
\cite{hrushovskimanuscript}) and of corollary 1.9 in
\cite{hrushovskimanuscript}.

Let $\mathcal{L}$ be the language of ordered fields, let
$\mathcal{L_-}$ be the language of fields and let $\mathcal{L_+}$
be the language of ordered fields with a predicate $F$
representing a pseudo real closed subfield of $D$. $T$ can then be
interpreted as a theory in the language $\mathcal{L_+}$. Notice
that for any model $M:=(D,F)\models T$, the theory 
$Th^{\mathcal{L}}(M)$ is 
the theory of real closed fields. From now on we only work with
bounded pseudo real closed fields and $T$ is understood to be the
theory of a bounded PRC-field.

For any element $a$ and any set $B$ in $D$, let $tp(a/B)$ and
$tp_+(a/B)$ be the type of $a$ over $B$ in the languages of
$\mathcal{L}$ and $\mathcal{L_+}$ respectively.

We define the \emph{field-definable closure} in the following way:
given some set $A$, the field definable closure of $A$
($dcl^{field}(A)$) is the set of elements that are definable using
quantifier free formulas in the language $\mathcal{L_-}(A)$. Note
that any field-definable element can be seen as the only element
of an irreducible variety, so $F$ is field-definably closed.

Let $M:=(D,F)$ be a model of the theory $T$ (with $F$ a bounded
PRC-field) and let $M(i):=(D(i),F(i))$. A substructure $M_0$ of
$M$ is said to be \emph{full} if it is algebraically closed and
$acl_{\mathcal{L_-}}(F^M)\subseteq dcl^{field}(F^M\cup
acl_{\mathcal{L_-}}(F^{M_0}))$.

\begin{remark}
If $M_1$ is a model of $T$ and $M_0$ is a full submodel of $M_1$,
then $M_0(i)$ is a full submodel of $M_1(i)$.
\end{remark}

\begin{proof}
We need to show that $M_0(i)$ is a submodel of $M_1(i)$, that
$D_j(i)$ is algebraically closed for $j=0,1$ and the ``fullness''
condition. Both $D_0$ and $D_1$ are real closed fields so both
$D_0(i)$ and $D_1(i)$ are algebraically closed.

The other conditions follow from the fact that all the field
operations inside $M_j(i)$ (and therefore the theory) are
interpretable inside $M_j$ using only quantifier free formulas in
$\mathcal{L_-}$.
\end{proof}

\begin{remark}
\label{1.4.1} If $M_1$ is a model of $T$ and $M_0$ is a submodel
of $M_1$, then $M_0$ is a full submodel of $M_1$.
\end{remark}

\begin{proof}
Let $M_0=(D_0, F^{M_0})$ be a submodel of $M_1$. By \ref{PRCPAC},
$M_0(i)$ is a submodel of $M_1(i)$ and by
\cite{hrushovskimanuscript} it is a full submodel so
\[ M_1(i)\subseteq dcl^{field}\left( F^{M_1}(i)\cup acl \left( F^{M_0}(i) \right) \right) \]
understanding $acl$ in the PAC-field sense. Once again, the
interpretability of the structure $F^{M_1}(i)\cup F^{M_0}(i)$ in
$F^{M_1}\cup F^{M_0}\restriction_{\mathcal{L_-}}$ proves the
remark.
\end{proof}

The following lemma is key for much of the rest of the proof. We
prove that the algebraic closure of a subset
of $F$ is in the field-definable closure of the union of $F$ and a
full submodel.

\begin{lemma}
\label{1.5} Let $M$ be a model of a (bounded PRC) theory $T$ and
$M_0$ a full submodel. Suppose $F^{M_0}\subseteq C\subseteq F^M$,
$acl(C)\cap F^M=C$ and $a\in acl(C)$ for some $a$. Then there
exists $e\in dcl^{field}(a,C)\cap M_0$ such that $a\in
dcl^{field}(e,C)$.
\end{lemma}
\begin{proof}
As with the previous two remarks we can prove the lemma 
using the analogue result for
PAC-fields (lemma 1.5 in \cite{hrushovskimanuscript}) and
interpretability of $C(i)$ and $M_0(i)$ by quantifier free
formulas in $\mathcal{L_-}(C)$ and $\mathcal{L_-}(M_0)$
respectively.
\end{proof}

\begin{proposition}
\label{1.6} Let $T,M$ and $M_0$ be as in the lemma above. Let $A$
be an algebraically closed subset of $M$ that contains $M_0$. Then
$T\cup diag(A)$ is complete (in the added language).
\end{proposition}

\begin{proof}
The proof is the same as the proof of proposition 1.6 in
\cite{hrushovskimanuscript} once we change stationary formulas to
(our corresponding) open subsets (in the order topology) of
irreducible affine varieties and $dcl$ to $dcl^{field}$.
\end{proof}

The following corollaries follow straight from the proofs of
corollaries 1.7, 1.8 and 1.9 in \cite{hrushovskimanuscript}.

\begin{corollary}
\label{1.7} $T$ is model complete.
\end{corollary}

\begin{corollary}
\label{1.8} A submodel $M_0$ of a model $M$ of $T$ is a full
submodel if and only if it is algebraically closed, full and
$F^{M_0}$ is a PRC subset of $M_0$.
\end{corollary}

\begin{corollary}
\label{1.9} Let $M$ be a model of $T$, $M_0$ a full submodel. Then
algebraic closure (in $\mathcal{L_+}$) over $M_0$ in $M$ and
coincides with field algebraic closure over $M_0$.
\end{corollary}

We can now prove the theorem.

\begin{proof}{(of theorem \ref{minimalityPRC})}

Let $M$ be a bounded large saturated PRC field, let $e\in
dcl^{eq}(F^M)$ be an imaginary element and let $N$ be an
elementary submodel of $F$. By lemma \ref{PRCPAC} we know that 
$N(i)$ is an
elementary submodel of $F(i)$. By \cite{hrushovskimanuscript}
corollary 3.2 (elimination of imaginaries for bounded PAC
structures) $e$ is interdefinable with some tuple $c\in F(i)$ in
the structure $(F(i),a)_{a\in N}$. By interpretability of $F(i)$
in $F\restriction \mathcal{L_-}$ we can find such a $c$ in $F$
which proves e.i. for $(F,a)_{a\in N}$.

Now, let $\overline{N}$ be the real closure of $N$. By corollary
\ref{1.8} $(N, \overline{N})$ is a full submodel of $M$ and by
corollary \ref{1.9} algebraic closure in $M$ and in
$M|\mathcal{L}$ coincide once we add constants for all the
elements in $N$. Therefore, algebraic closure in $(F,a)_{a\in N}$
is the field algebraic closure so it satisfies Steinitz exchange
property.

By theorem \ref{uth1}, $T^*$ is \th-minimal. But being \th-minimal
is a property that is invariant under adding constants, so $T$ is
\th-minimal.
\end{proof}

\subsection{Model Companions of Large differential Fields of
Characteristic 0}

In \cite{tressl}, Marcus Tressl introduced a first order theory in
the language of differential rings with $k$ derivatives called
$UC$\footnote{$UC$ is a system of axioms that basically says that
any algebraically prepared system (a system of differential
polynomials that is consistent with $T^{field}$) has a
realization.} with the property that if $T$ is the theory of a
differential field (in the language of differential rings) such
that the restriction $T^{field}$ of $T$ to the language of fields
is model complete and has large models, then $T^*:=T^{field}\cup
UC$ is the model companion of $T$. We prove in this section that if
$T^{field}$ is rosy, $T^*$ is rosy. As a corollary, we prove that
if $T^{field}$ is stable then $T^*$ is stable.

In this section we work with a theory $T$ in the language of
differential rings containing the theory of differential fields
such that $T^{field}$ is a model complete theory with large
models. We assume the reader has familiarity with the results and
definitions in \cite{tressl}.

\begin{definition}
Given a complete type $p(x)$, let $p^{field}(x,x_1,x_2, \dots ,
x_k x_{11}, \dots)$ be the restriction of $p(x)$ to the language
of rings obtained by replacing all the derivatives of $x$ by free
variables and let $p^{qf}(x)$ be the set of all quantifier free
formulas in $p(x)$.
\end{definition}

\begin{lemma}
\label{tressl} Let $p(x)$ be a complete type over some set $A$.
Then, if $p^{field}\left(\bar{x}\right)$ is consistent
with $T^{field}$ and $p^{qf}(x)$ is a consistent differential
type, there is some model $N\models T^*$ containing $A$ such that
$p(x)$ is realized in $N$.
\end{lemma}

\begin{proof}
We can assume without loss of generality that $A$ is a
differential field.

By \cite{tressl} theorem $7.1 (ii)$, $p(x)$ is implied by
$p^{field}\left(\bar{x}\right)\cup p^{qf}(x)$. Let
$A^c$ be the differential closure of $A$ (see \cite{mcgrail}) and
let $M^{field}$ be a model of $T^{field}$ such that there is some
tuple $\left<a, a_1, a_2, \dots , a_k, a_{11} \dots \right>$ in
$M^{field}$ realizing $p^{field}\left( \bar{x}\right)$.
Let $F$ be the subfield of $M$ generated by $A\cup \left\{
\bar{a_I}\right\}$.

Since $A^c$ is differentially closed there is some $a'\in M^c$
realizing $p^{qf}(x)$. Let $D$ be the differential subfield of
$M^c$ generated by $\{Aa'\}$. By definition $a'$ satisfies
$p^{qf}(x)$ so the tuple $\left<a',d_1a',d_2a', \dots , d_k a',
d_1d_1a' \dots\right>$ satisfies all the quantifier free formulas
that appear in $p^{field}$ and we have a map of rings sending
$\left< A, a', d_1a', d_2a', \dots , d_k a', d_1d_1a' \dots
\right>$ to $\left< A, a, a_1, a_2, \dots, a_k,  a_{11} \dots
\right> $. We can use this map to equip $F$ with $k$ commuting
derivatives extending those in $A$ in such way that $a\models
p^{qf}(x)$. By \cite{kolchin} we can define $k$ commuting
derivatives in $M^{field}$ extending those in $F$ so there is 
some differential field $M$ such that $F$ is a subfield of
$M$ and such that $M^{field}$ is the restriction of $M$ to the
language of rings. By \cite{tressl} theorem 6.2(II) there is a
differential field $L$ extending $M$ satisfying $T^*$. By model
completeness of $T^{field}$, the field type of $a$ over $A$ is the
same in $M$ as it is in $L$ so $L\models p^{field}(a)$.

Therefore, $L\models p^{field}(a)\cup p^{qf}(a)$ and $L\models
p(a)$.
\end{proof}

\begin{corollary}
\label{strdividing} Let $T, T^*, T^{field}$ be as above, $M$ be a
model of $T^*$, $A\subset B$ subsets of $M$. Let $p(x,B)$ be a
type over $B$ realized by some $a$ and let $\text{k-DCF}$ be the
theory of closed differential fields with $k$ commuting
derivatives. Then, if $T^*$ implies that $p(x)$ strongly divides
over $A$, either $T^{field}$ implies $p^{field}$ strongly divides
over $A$ or $\text{k-DCF}$ implies $p^{qf}$ strongly divides over
$A$\footnote{If we look in the differential closure of $M$ the
quantifier free type of $a$ over $B$ is precisely $p^{qf}$; by
$p^{qf}$ strongly dividing over $A$ in the sense of
$\text{\emph{k}-DCF}$ we mean that in the differential closure of
$M$ the quantifier free type of $a$ over $B$ strongly divides over
$A$. This is equivalent to say that the differential ideal of
$p\qf$ has smaller differential rank than the differential ideal
of $p\qf \mid _A$ (see remark below).}.
\end{corollary}

\begin{proof}
Suppose $p(x)$ strongly divides over $A$ in the sense of $T^*$ and
let $q(Y,A)$ be the type of $B$ over $A$. By definition, there is
some $n\in \mathbb{N}$ such that
\[ \bigcup_{i=1}^n p(x,Y_i) \cup \bigcup_{i=1}^n q(Y_i,A) \] is
inconsistent with $T^*$ and $q(Y,A)$ is non algebraic. By lemma
\ref{tressl}, either
\[ \bigcup_{i=1}^n p^{field}(x,Y_i) \cup \bigcup_{i=1}^n q^{field}(Y_i,A)
\]is inconsistent with $T^{field}$, or
\[ \bigcup_{i=1}^n p^{qf}(x,Y_i) \cup \bigcup_{i=1}^n q^{qf}(Y_i,A)
\] is inconsistent with $\text{\emph{k}-DCF}$. Since on one hand
$tp^{field}(B/A)$ is precisely $q^{field}(Y,A)$ and on the other
$q^{qf}(Y,A)$ is quantifier free and we have elimination of
quantifiers in $\text{\emph{k}-DCF}$, the result follows.
\end{proof}

\begin{remark}
Note that strong dividing implies forking so, by \cite{mcgrail},
if $\text{k-DCF}$ implies $p^{qf}$ strongly divides over $A$ then
the $\Delta$-differential rank of the prime ideal $\mathcal{I}_p$
generated by differential polynomials in $p$ is smaller than the
differential rank of the prime ideal $\mathcal{I}_{p\mid A}$
generated by the differential polynomials in $p$ with coefficients
in $A$.
\end{remark}

\begin{theorem}
If $T^{field}$ is rosy so is $T^*$.
\end{theorem}

\begin{proof}
Let $\mathcal{M}$ be a monster model of $T^*$, let $A\subset B$ be
(small) subsets of $\mathcal{M}$ and let $p(x)$ be a type over $B$
which \th-forks over $A$. By definition of \th-forking this is
witnessed by
\[p(x)\vdash \bigvee_{i=1}^n \phi_i(x,b_i)\] where for each of $i$ we have
that $\phi_i(x,b_i)$ strongly divides over $Ac_i$ for some $c_i$;
let $D:=\bar{b_i}$.

\begin{claim}
Let $p\field(x)$ and $p\qf(x)$ be as above. Then either
$p\field(x)$ \th-forks over $A$ in the sense of $T\field$ or
$p\qf(x)$ \th-forks over $A$ in the sense of $\text{k-DCF}$.
\end{claim}

\begin{proof}
Let $p\field(x,D)$ and $p\qf(x,D)$ be types over $B\cup D$ which 
are, respectively, non \th-forking extensions of $p\field(x)$ and
$p\qf(x)$ in the sense of $T\field$ and $k$-DCF. By lemma
\ref{tressl} the type $p\field(x,D) \cup p\qf(x,D)$ is consistent
with $T^*$ so it is realized by some $a\in M$. By \cite{tressl}
theorem 7.1(ii) we know that $p\field(x,D) \cup p\qf(x,D)$ implies
$tp(a/BD)$ so by construction there is some $i$ such that $a$
realizes $\phi_i(x,b_i)$. This implies that $tp(a/BD)$ strongly
divides over $Ac_i$. By theorem \ref{strdividing}, either
$p\field(x,D)$ \th-divides over $A$ in the sense of $T\field$ or
$p\qf(x,D)$ \th-divides over $A$ in the sense of
$k$-DCF. By construction neither of them \th-forks
over $B$ so transitivity of \th-forking implies that either
$p\field(x)$ \th-forks over $A$ in the sense if $T\field$ or
$p\qf(x)$ \th-forks over $A$ in the sense of differentially closed
fields.
\end{proof}

To finish the prove of the theorem, just note that if we had an
infinite \th-forking chain in a model of $T^*$ we would have an
infinite \th-forking chain in a model of $T\field$ or some
differential ideal with infinite differential rank. Since
$T\field$ is rosy and $\text{k-DCF}$ is stable, this cannot
happen.
\end{proof}

\begin{corollary}
If $T\field$ is stable so is $T^*$.
\end{corollary}

\begin{proof}
Let $M$ be a monster model of $T^*$. Since \th-forking is an
independence relation it is enough to show that given a small
model $N\models T^*$ a type $p(x)$ over $N$ and some tuple $a$
there is a unique non \th-forking extension of $p(x)$ to $Na$. Let
$q(x,a)$ and $r(x,a)$ be two non \th-forking extensions of $p(x)$.
Let $q\field, r\field, p\field, q\qf, r\qf$ and $p\qf$ be the
types obtained by restricting types $q,r$ and $p$ to the language
of rings and to the quantifier free formulas. By quantifier
elimination and stability of $k$-DCF we have that
$r\qf=q\qf$ and by stability of $T\field$ (and therefore
uniqueness of non \th-forking extensions) $r\field=q\field$. Since
both $q(x,a)$ and $r(x,a)$ are implied by $q\field \cup q\qf$ and
 $r\field \cup r\qf$ respectively, we conclude that
$r(x,a)=q(x,a)$ so that there is a unique non \th-forking extension of
$p(x)$ to $Na$.
\end{proof}

\nocite{vandenDries, Kim/Pillay/stableforking}

\bibliographystyle{alpha}

\bibliography{arti2}

\begin{thebibliography}{GvdDM81}

\bibitem[Cha99]{chatzidakis}
Zo{\'e} Chatzidakis.
\newblock Simplicity and independence for pseudo-algebraically closed fields.
\newblock In {\em Models and computability (Leeds, 1997)}, volume 259 of {\em
  London Math. Soc. Lecture Note Ser.}, pages 41--61. Cambridge Univ. Press,
  Cambridge, 1999.

\bibitem[CP98]{chatzidakis/pillay}
Z.~Chatzidakis and A.~Pillay.
\newblock Generic structures and simple theories.
\newblock {\em Ann. Pure Appl. Logic}, 95(1-3):71--92, 1998.

\bibitem[GvdDM81]{cherlin/van/mc}
Gregory Gherlin, Lou van~den Dries, and Angus Macintyre.
\newblock Decidability and undecidability theorems for pac-fields.
\newblock {\em Bull. Amer. Math. Soc.}, 4(1):101--104, 1981.

\bibitem[HKP00]{Hart/Kim/Pillay}
Bradd Hart, Byunghan Kim, and Anand Pillay.
\newblock Coordinatisation and canonical bases in simple theories.
\newblock {\em J. Symbolic Logic}, 65(1):293--309, 2000.

\bibitem[Hru91]{hrushovskimanuscript}
Ehud Hrushvski.
\newblock Pseudo-finite fields and related structures.
\newblock {\em manuscript}, 1991.

\bibitem[Kim96]{Kim}
Byunghan Kim.
\newblock {\em Simple First Order Theories}.
\newblock PhD thesis, University of Notre Damme, 1996.

\bibitem[Kol73]{kolchin}
E.R. Kolchin.
\newblock Differential algebra and algebraic groups.
\newblock {\em Pure and Applied Mathematics}, 54, 1973.

\bibitem[KP01]{Kim/Pillay/stableforking}
Byunghan Kim and A.~Pillay.
\newblock Around stable forking.
\newblock {\em Fund. Math.}, 170(1-2):107--118, 2001.
\newblock Dedicated to the memory of Jerzy \L o\'s.

\bibitem[McG00]{mcgrail}
Tracey McGrail.
\newblock The model theory of differential fields with finitely many commuting
  derivatives.
\newblock {\em J. Symbolic Logic}, 65(2):885--913, 2000.

\bibitem[Ons]{Onshuus}
Alf Onshuus.
\newblock Properties and consequences of th-forking.
\newblock {\em Preprint}.

\bibitem[Ons02]{phd}
Alf Onshuus.
\newblock {\em Thorn-Forking in Rosy Theories}.
\newblock PhD thesis, University of California at Berkeley, 2002.

\bibitem[Pre81]{Prestel}
Alexander Prestel.
\newblock Pseudo real closed fields.
\newblock In {\em Set theory and model theory}, volume 872 of {\em Lecture
  Notes in Math.}, pages 127--156. Springer, Berlin, 1981.

\bibitem[Sch00]{Schmid}
Joachim Schmid.
\newblock A density property for pseudo real closed fields.
\newblock {\em Konstanzer Schriften in Mathematik und Informatik}, (124), 2000.

\bibitem[Tre]{tressl}
Marcus Tressl.
\newblock A uniform companion for large differential fields of characteristic
  0.
\newblock {\em Preprint}.

\bibitem[vdD98]{vandenDries}
Lou van~den Dries.
\newblock {\em Tame topology and o-minimal structures}, volume 248 of {\em
  London Mathematical Society Lecture Note Series}.
\newblock Cambridge University Press, Cambridge, 1998.

\bibitem[Wag00]{Wagner}
Frank~O. Wagner.
\newblock {\em Simple theories}, volume 503 of {\em Mathematics and its
  Applications}.
\newblock Kluwer Academic Publishers, Dordrecht, 2000.

\end{thebibliography}


\end{document}